# CHAIN DECOMPOSITION THEOREMS FOR ORDERED SETS
AND OTHER MUSINGS

JONATHAN DAVID FARLEY

*This paper is dedicated to the memory of Prof. Garrett Birkhoff*

ABSTRACT. A brief introduction to the theory of ordered sets and lattice theory is given. To illustrate proof techniques in the theory of ordered sets, a generalization of a conjecture of Daykin and Daykin, concerning the structure of posets that can be partitioned into chains in a "strong" way, is proved. The result is motivated by a conjecture of Graham's concerning probability correlation inequalities for linear extensions of finite posets.

## 1. INTRODUCTION

Order has played a rôle in human civilization for as long as the North Star has hung above our heads. The theory of ordered sets, however, is a relatively new discipline.

Lattice theory and the theory of ordered sets are exciting areas with a number of surprising connections with other branches of mathematics, including algebraic topology, differential equations, group theory, commutative algebra, graph theory, logic, and universal algebra. Both fields have many important applications, for example, to scheduling problems, the semantics of programming languages, the logic of quantum mechanics, mathematical morphology and image analysis, circuit design, and cryptography [2], [10], [14], [21].

Below we present a sampling of theorems — some famous, some not — dealing with decompositions of ordered sets into chains. We will expand upon our own original work (presented in §7) in a future note.

1991 *Mathematics Subject Classification.* Primary 06A07, 05A18; Secondary 06-02, 06-06, 05-02, 05-06.

The author would like to thank Dr. William Massey and Dr. Nathaniel Dean for inviting him to speak at the Second Conference for African-American Researchers in the Mathematical Sciences. The author would also like to thank Prof. Bill Sands for permitting him to publish his theorem. This paper is based on a lecture given at the Institute for Advanced Study in Princeton, New Jersey, on June 27, 1996.

Typeset by $\mathcal{A}_{\mathcal{M}}\mathcal{S}$-TEX

2. Basic terminology, notation, and examples

A *partially ordered set*, or *poset*, is a set with a binary relation $\leqslant$ that is reflexive, transitive, and antisymmetric. Elements $a, b \in P$ are *comparable*, denoted $a \sim b$, if $a \leqslant b$ or $a \geqslant b$; otherwise they are *incomparable*, denoted $a \parallel b$. An element $a$ is a *lower cover* of $b$ if $a < b$ and $a < c \leqslant b$ implies $b = c$ for all $c \in P$; in this case, $b$ is an *upper cover* of $a$, which we denote by $a \prec b$. The *Hasse diagram* of a finite poset is a graph in which each vertex represents an element, and an edge drawn upward from $a$ to $b$ means that $a \prec b$. Hence $a \leqslant b$ if and only if one can go from $a$ to $b$ by tracing the edges upward.

A *chain* is a subset $C$ of $P$ such that $a \sim b$ for all $a, b \in C$. A chain $c_1 < c_2 < \cdots < c_n$ is *saturated* if $c_1 \prec c_2 \prec \cdots \prec c_n$. The *length* of a finite chain is $\#C - 1$ (where $\#C$ is the cardinality of $C$). The *height* $\operatorname{ht} x$ of an element $x$ is the supremum of the lengths of the finite chains with greatest element $x$. The height $\operatorname{ht} P$ of a poset is the supremum of the heights of its elements. A finite poset is *ranked* if, for all $x \in P$, every chain with greatest element $x$ is contained in a chain of length $\operatorname{ht} x$ with greatest element $x$. A saturated chain $c_1 \prec c_2 \prec \cdots \prec c_n$ in a ranked poset is *symmetric* if $\operatorname{ht} c_1 = \operatorname{ht} P - \operatorname{ht} c_n$.

An *antichain* is a subset $A$ of $P$ such that $a \parallel b$ for distinct $a, b \in A$. The *width* of a poset is the supremum of the cardinalities of its antichains.

For all $p \in P$, let $\uparrow p := \{ q \in P \mid p \leqslant q \}$ and $\downarrow p := \{ q \in P \mid p \geqslant q \}$. An element $p$ is *maximal* if $\uparrow p = \{p\}$; let $\operatorname{Max} P$ denote the set of maximal elements. An *up-set* of $P$ is a subset $U$ such that, for all $u \in U$, $\uparrow u \subseteq U$.

The *disjoint sum* of two posets $P$ and $Q$ is the poset with underlying set $P \cup Q$ and the inherited order on $P$ and $Q$, but no comparabilities between elements of $P$ and elements of $Q$. The *ordinal sum* of $P$ and $Q$ is the poset with underlying set $P \cup Q$ and the inherited order on $P$ and $Q$, but with $p < q$ for all $p \in P$ and $q \in Q$. It is denoted $P \oplus Q$.

A *lattice* is a non-empty poset $L$ such that, for all $a, b \in L$, the least upper bound of $a$ and $b$ exists, called the *join* of $a$ and $b$ (denoted $a \vee b$), and the greatest lower bound of $a$ and $b$ exists, called the *meet* of $a$ and $b$ (denoted $a \wedge b$). A *sublattice* is a non-empty subset closed under join and meet. A lattice is *distributive* if, for all $a, b, c \in L$, $a \wedge (b \vee c) = (a \wedge b) \vee (a \wedge c)$ [equivalently, for all $a, b, c \in L$, $a \vee (b \wedge c) = (a \vee b) \wedge (a \vee c)$].

An example of a chain is the real line with the usual ordering. The power set of a set is a distributive lattice; here join corresponds to set union and meet to set intersection. The collection of all subsets with 42 elements is an antichain. The collection of subsets of $\{x, y, z\}$ containing either $\{x, y\}$ or $\{z\}$ is an up-set of the power set of $\{x, y, z\}$.

The left and right sides of Figure 1 are the Hasse diagrams of the four-element fence $P$ and the four-element crown $Q$, respectively. Their disjoint sum and ordinal sum are shown in Figures 1 and 2.

Another example of a distributive lattice is the set of natural numbers ordered by divisibility ($a \leqslant b$ if $a$ divides $b$). In this lattice, $a \vee b$ is the least common multiple of $a$ and $b$, and $a \wedge b$ is the greatest common divisor of $a$ and $b$.

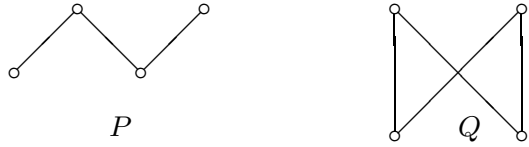

**Figure 1.** The disjoint sum of $P$ and $Q$.

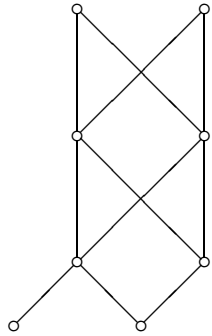

**Figure 2.** The ordinal sum $P \oplus Q$.

The lattice of equivalence classes of propositions in propositional logic is a distributive lattice when ordered by implication (the class $[p]$ of the proposition $p$ is less than or equal to $[q]$ if $p$ implies $q$). In this lattice, $[p] \vee [q] = [p \text{ or } q]$ and $[p] \wedge [q] = [p \text{ and } q]$.

The lattice of subspaces of a vector space is *not* distributive if the dimension of the space is at least 2. The five-element lattices $M_3$ and $N_5$ of Figure 3 are also non-distributive.

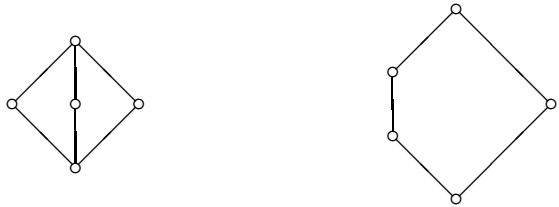

**Figure 3.** The lattices $M_3$ and $N_5$.

Figure 4 shows the Hasse diagram of a poset that is not a lattice.

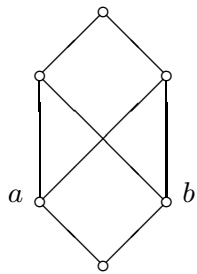

**Figure 4.** A non-lattice.

The elements $a$ and $b$ have common upper bounds, but not a *least* upper bound, so the poset is not a lattice.

Basic references on lattice theory and the theory of ordered sets are [1] and [4].

## 3. Dilworth's theorem

Even mathematicians who know very little about the theory of ordered sets have heard of Dilworth's Theorem. (R. P. Dilworth actually proved many beautiful and important theorems in the theory of ordered sets, all of them called "Dilworth's Theorem." We, in particular, should note that Dilworth helped to develop mathematics programs on the Continent [3].)

Suppose one desires to partition a poset into chains, using the smallest number possible. (For instance, one would want such a partition for some scheduling applications.) If $P$ has finite width $w$, then clearly we cannot do better than $w$ chains. Dilworth's Theorem asserts that we can achieve this bound ([8], Theorem 1.1).

**Dilworth's Theorem.** *Let $P$ be a poset of finite width $w$. Then there exists a partition of $P$ into $w$ chains, and this is best possible.*

Note that we are not assuming $P$ is finite; Dilworth's Theorem fails, however, for posets of infinite width, even if every antichain is finite [17].

For example, in Figure 5, the poset has width 4 ($\{a, b, e, f\}$ is an antichain), and can be partitioned into 4 chains (e.g., $\{0, a, d, 1\}$, $\{b\}$, $\{e\}$, $\{c, f\}$).

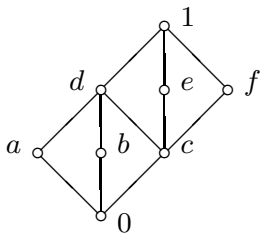

**Figure 5.** An illustration of Dilworth's Theorem.

For more on Dilworth's work, see [9].

## 4. Birkhoff's theorem

Garrett Birkhoff was one of the pioneers of lattice theory and universal algebra. Two of his many significant contributions involve distributive lattices.

Birkhoff characterized distributive lattices in terms of forbidden sublattices [analogous to Kuratowski's characterization of planar graphs in terms of forbidden subgraphs ([11], Theorem 6.2.1)].

**Theorem** (Birkhoff–Dedekind). *A lattice is distributive if and only if neither $M_3$ nor $N_5$ is a sublattice.*

The theorem makes it easy to spot non-distributive lattices. For instance, the lattice of Figure 5 is not distributive. For a proof of the theorem, see [4], 6.10.

Now we know what lattices are *not* distributive. Which lattices *are* distributive?

Here is one way to construct distributive lattices: Take a finite poset, and order its up-sets by inclusion. This procedure yields a finite distributive lattice, in which join and meet are given by union and intersection, respectively.

For example, if Abubakari is the three-element fence (Figure 6), its lattice of up-sets is shown in Figure 7.

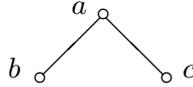

**Figure 6.** The poset Abubakari.

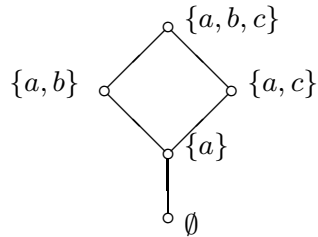

**Figure 7.** The distributive lattice of up-sets of Abubakari.

Another theorem of Birkhoff's asserts that *all* finite distributive lattices arise in this way ([4], 8.17). One consequence, which we shall need later, is that finite distributive lattices are ranked.

**Theorem** (Birkhoff). *Every finite distributive lattice is the lattice of up-sets of a finite poset.*

There is an analogous theorem for infinite distributive lattices, but, instead of finite posets, one must use compact Hausdorff partially ordered spaces, called *Priestley spaces* [18].

## 5. Sands' matching theorem

Now we may prove an amusing theorem due to Bill Sands [20]:

**Theorem** (Sands, 197∗). *Every finite distributive lattice with an even number of elements can be partitioned into two-element chains.*

In the terminology of graph theory, the theorem asserts that the *comparability graph* of a finite distributive lattice with an even number of elements [in which $(a,b)$ is an edge if $a \sim b$] has a perfect matching without loops ([11], §7.1).

The reader might like to attempt his or her own proof before proceeding.

*Proof.* We prove the theorem by induction on the size of $L$. If $\#L$ is odd, and $L$ has greatest element 1, we prove that $L \setminus \{1\}$ can be partitioned into two-element chains.

Represent the lattice as the lattice of up-sets of a finite poset $P$. Pick $p \in P$, and consider the the elements $a := \uparrow p$ and $b := P \setminus \downarrow p$ of $L$. Then $L$ is partitioned into the finite distributive lattices $\uparrow a$ and $\downarrow b$. There are four cases to consider, depending on the parities of $\# \uparrow a$ and $\# \downarrow b$. □

The theorem is interesting because most chain-decomposition theorems apply either to arbitrary posets (e.g., Dilworth's Theorem) or only to power set lattices (e.g., [16]). Sands' Matching Theorem is obvious for power set lattices, but does not apply to all posets, not even all lattices. For example, $M_4$ (Figure 8) is non-distributive, by Birkhoff's theorem; it cannot be partitioned into two-element chains. (This example is due to Sands.)

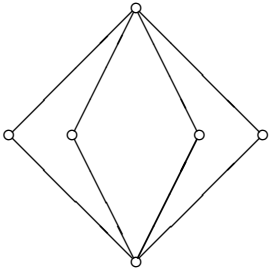

**Figure 8.** The lattice $M_4$ cannot be partitioned into two-element chains.

## 6. The parable of the tennis players

Our main theorem is a generalization of a conjecture of Daykin and Daykin that deals with special sorts of chain-partitions. The motivation for the conjecture comes from a parable, due to Graham, Yao, and Yao ([13], §1):

Imagine there are two teams of tennis players, $A$ and $B$. The players of Team $A$ are linearly ranked from best to worst, as are the players from Team $B$, but we know only in a few cases how individual players from Team $A$ compare with individuals from Team $B$.

We may describe this set-up using a poset: The players are the elements, and a relation $p < q$ means that player $p$ is worse than player $q$. Hence the subset corresponding to Team $A$ is a chain, as is the subset corresponding to Team $B$ (so that we have a poset of width at most 2). For an example, see Figure 9.

In Figure 9, we do not know whether Asmodeus (in Team $A$) is better or worse than Beezelbub (in Team $B$). We might ask for the *probability* that Asmodeus is worse than Beezelbub, $Pr$(Asmodeus<Beezelbub). We calculate this probability by looking at all the possible linear rankings of the players in both teams that are consistent with what we already know about which players are better than which.

Formally, we are looking at all the possible *linear extensions* of the poset, the bijective order-preserving maps from $P$ into a chain of cardinality $\#P$. By counting the number of these in which (the image of) Asmodeus is below (the image of) Beezelbub, and dividing by the total number of linear extensions (assuming all are equally likely), we have $Pr$(Asmodeus<Beezelbub).

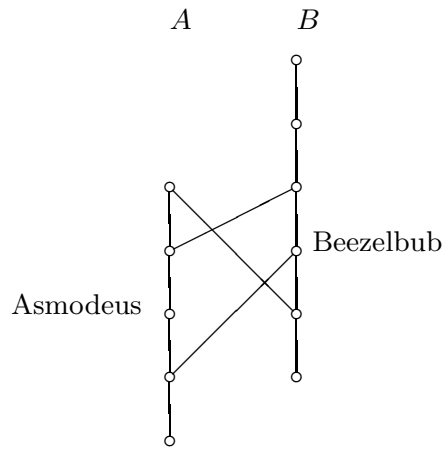

**Figure 9.** A tale of two tennis teams.

Now suppose we are given additional information. Namely, suppose we learn that certain players from Team $A$ are worse than certain players from Team $B$. (Perhaps the two teams have just finished playing a tournament.) This information supports the idea that the players from Team $A$ are worse than the players from Team $B$, making it more likely that Asmodeus is, in fact, worse than Beezelbub.

Formally, we would expect:

$$Pr(\text{Asmodeus}<\text{Beezelbub}) \leqslant Pr(\text{Asmodeus}<\text{Beezelbub} \mid a < b \text{ \& } \cdots \text{\& } a' < b'),$$

the conditional probability that Asmodeus is worse than Beezelbub given that $a$ is worse than $b$, $a'$ worse than $b'$, etc. (See Figure 10.)

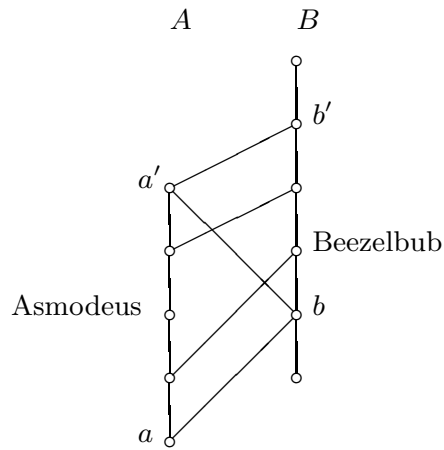

**Figure 10.** The poset of tennis players after a tournament.

We can formalize a more general situation. Let $P$ be a poset partitioned into subsets $A$ and $B$ (not necessarily chains). Suppose that, whenever $x$ and $y$ are disjunctions of statements of the form

$$a < b \text{ \& } \cdots \text{\& } a' < b'$$

we know that
$$Pr(x)Pr(y) \leqslant Pr(x\&y).$$

Then we say the partition has the *positive correlation property*.

Using a result of Daykin and Daykin ([5], Theorem 9.1), a question of Graham's ([12], pp. 232-233) becomes:

**Question** (Graham). *Let $P = A \cup B$ be a partition of a finite poset such that, for all $a \in A$ and $b \in B$,*
$$a < b \text{ implies } P = \uparrow a \cup \downarrow b$$
*and*
$$a > b \text{ implies } P = \downarrow a \cup \uparrow b.$$
*Does the partition have the positive correlation property?*

## 7. A CONJECTURE OF DAYKIN AND DAYKIN

In light of Graham's question, Daykin and Daykin ([5], §9) conjectured the following:

**Conjecture** (Daykin–Daykin). *Let $P$ be a finite poset partitioned into chains $T_1$, $T_2$, and $T_3$ such that, if $p$ and $q$ are in different chains and $p < q$, then $P = \uparrow p \cup \downarrow q$.*

*Then $P$ is an ordinal sum $R_1 \oplus \cdots \oplus R_n$ such that, for $1 \leqslant i \leqslant n$, either*

(1) *$R_i$ is disjoint from some $T_j$, $j \in \{1, 2, 3\}$, or*
(2) *for all $p$, $q \in R_i$, if $p$ and $q$ are in different chains, then $p \parallel q$.*

With this conjecture, Daykin and Daykin asserted that certain correlation inequalities would follow.

Figure 11 shows an example of a poset satisfying the hypotheses of the conjecture. Let $T_1$ be the chain of all the elements on the left, $T_2$ all the elements on the right, and $T_3$ the remainder.

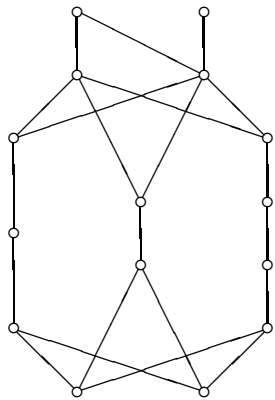

**Figure 11.** A poset satisfying the hypotheses of the Daykin–Daykin conjecture.

Note that the conclusion of the conjecture implies that $P$ is the ordinal sum of width 2 posets and disjoint sums of chains.

Tseng and Horng [23] claim to have a proof, but their proof rests on a mistaken assertion. The author has proven a more general statement, with a much simpler proof, which applies to infinite posets; we shall write another note on this theorem. (Incidentally, although the fact was not publicized, it is reported in [6], p. 84 that J. M. Robson proved the Daykin–Daykin conjecture by an inductive argument.)

The results in the sequel are due to the author.

By $P(R_1, \ldots, R_n; T_1, \ldots, T_s)$, we mean:

(1) $P$ is a finite poset;
(2) $P = R_1 \oplus \cdots \oplus R_n$;
(3) $P = T_1 \cup \cdots \cup T_s$;
(4) $\{T_j \mid j = 1, \ldots, s\}$ are disjoint chains;
(5) for $i \in \{1, \ldots, n\}$, either
$R_i$ is disjoint from $T_j$ for some $j \in \{1, \ldots, s\}$ or
$p \parallel q$ for all $p \in R_i \cap T_j$ and $q \in R_i \cap T_k$ ($1 \leqslant j < k \leqslant s$).

**Lemma.** *Let $P'(R'_1, \ldots, R'_m; T'_1, T_2, \ldots, T_s)$ where $s \geqslant 3$. Suppose that $P = P' \cup \{x\}$ is a poset with a new element $x$, $\operatorname{ht} x = \operatorname{ht} P$ and for all $p \in P$ such that $\operatorname{ht} p = \operatorname{ht} P$,*

$$\#\{q \in P \mid q \prec p\} \leqslant \#\{q \in P \mid q \prec x\}.$$

*Assume that $T_1 := T'_1 \cup \{x\}$ is a chain. Finally, assume that, for all $p \in T_j$, $q \in T_k$, $1 \leqslant j, k \leqslant s$, $j \neq k$, $p < q$ implies $P = \uparrow p \cup \downarrow q$.*

*Then $P(R_1, \ldots, R_n; T_1, \ldots, T_s)$ for some $R_1, \ldots, R_n \subseteq P$.*

*Proof.* If $\operatorname{Max} P = \{x\}$, let $R_i := R'_i$ ($1 \leqslant i \leqslant m$) and $R_{m+1} := \{x\}$. We may thus assume that $\operatorname{Max} P \neq \{x\}$.

We claim that $p \leqslant x$ for all $p \in R'_{m-1}$.

For if $p \not\leqslant x$, where $p \in R'_{m-1} \cap T_k$, then $\uparrow p \subseteq T_k$; in particular, $R'_m \subseteq T_k$. If $C \subseteq P$ is a chain containing $x$, then $C \cup \{y\} \setminus \{x\}$ is a chain for all $y \in R'_m$, so that, for some $y \in R'_m$, $\operatorname{ht} y = \operatorname{ht} P$. Hence $y$ has more lower covers than $x$, a contradiction.

If $R'_m \cup \{x\}$ is disjoint from some $T_j$, we finish by letting $R_i := R'_i$ ($1 \leqslant i < m$) and $R_m := R'_m \cup \{x\}$.

We may therefore assume that $(R'_m \cup \{x\}) \cap T_j \neq \emptyset$ for $j = 1, \ldots, s$.

*Case 1.* For $1 \leqslant j \leqslant s$, $R'_m \cap T_j \neq \emptyset$.

Assume there exists $p \in T_2 \cap R'_m$ such that $p < x$. Then $R'_m \subseteq \uparrow p \cup \downarrow x$ and $\uparrow p \subseteq T_2 \cup \{x\}$, so that $(\operatorname{Max} P) \setminus \{x\} \subseteq T_2$. Hence $R'_m \cap T_3 \cap \downarrow x \neq \emptyset$, so that $(\operatorname{Max} P) \setminus \{x\} \subseteq T_3$. Hence $\operatorname{Max} P = \{x\}$.

*Case 2.* For some $j \in \{1, \ldots, s\}$, $R'_m \cap T_j = \emptyset$.

Clearly $j = 1$. If $x$ has no lower covers in $R'_m$, then $R'_m$ is an antichain and we are done. Otherwise, let $p \in R'_m$ be such that $p \prec x$ and $\{p, x\}$ belongs to a chain of length ht $P$. Then for all $y \in \mathrm{Max}\, P$, $p \prec y$, and every lower cover of $y$ is below $x$.

Now let $q \in R'_m$ be a lower cover of $x$. For all $y \in \mathrm{Max}\, P$, there exists $z \in P$ such that $q \prec z \leqslant y$ and $p < z$, so that $z = y$. Therefore we may let $R_i := R'_i$ ($1 \leqslant i < m$), $R_m := R'_m \setminus \mathrm{Max}\, P$, and $R_{m+1} := \mathrm{Max}\, P$. □

**Theorem.** *Let $P$ be a finite poset and $T_1, \ldots, T_s$ disjoint chains covering $P$ ($s \geqslant 3$). Assume that, for all $p \in T_j$, $q \in T_k$ such that $p < q$, $P = {\uparrow}p \cup {\downarrow}q$ whenever $1 \leqslant j, k \leqslant s$ and $j \neq k$.*

*Then there exist subsets $R_1, \ldots, R_n$ of $P$ such that $P = R_1 \oplus \cdots \oplus R_n$ and, for $i = 1, \ldots, n$, either*

(1) $R_i \cap T_j = \emptyset$ *for some* $j \in \{1, \ldots, s\}$, *or*
(2) $p \parallel q$ *whenever* $p \in R_i \cap T_j$, $q \in R_i \cap T_k$, *and* $j \neq k$.

*Proof.* The theorem follows from the lemma by induction. □

The theorem will be generalized in a future note.

## 8. An open problem regarding chain decompositions of distributive lattices

Dilworth's Theorem asserts that, for every width $w$ poset, there *exists* a partition into $w$ chains. Both practical and philosophically-minded persons might prefer an explicit construction for such a minimal partition.

It is clear that any partition of a ranked poset into symmetric chains is a minimal chain-partition. Such a decomposition also tells us that the width of the poset is the number of elements of height $\lfloor \frac{1}{2} \mathrm{ht}\, P \rfloor$.

Examples of distributive lattices admitting such symmetric chain decompositions are finite power set lattices and lattices of divisors of positive integers [7].

Let $L(m, n)$ denote the poset of $m$-tuples $(x_1, \ldots, x_m)$ where $0 \leqslant x_1 \leqslant \cdots \leqslant x_m \leqslant n$, ordered as follows: $(x_1, \ldots, x_m) \leqslant (y_1, \ldots, y_m)$ if $x_i \leqslant y_i$ for $i = 1, \ldots, m$. This poset is a distributive lattice. Figure 12 shows $L(2, 3)$.

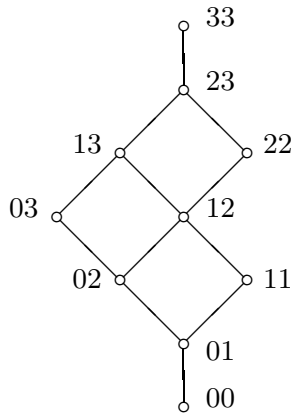

**Figure 12.** The lattice $L(2, 3)$.

**Conjecture** (Stanley). *The lattice $L(m, n)$ admits a symmetric chain decomposition.*

It is known that the conjecture holds for $L(3, n)$ and $L(4, n)$ [15], [19], [24]; it is also known that the width of $L(m, n)$ is indeed the number of elements of height $\lfloor \frac{1}{2} \operatorname{ht} L(m, n) \rfloor$, but the proof relies on heavy machinery from algebraic geometry ([22], §7). It would be interesting to have an elementary proof even of this fact.

## 9. Afterglow

"The Universe," said the philosopher with a glint in his eye, "is a big place." One might make a similar remark regarding the theory of ordered sets. What we have done is survey a tiny portion of the night sky. We hope that we have instilled in you, the reader, the desire to explore the constellations of lattice theory and the theory of ordered sets. Just perhaps, you may discover a star of your own.

Mathematical Sciences Research Institute, 1000 Centennial Drive, Berkeley, California 94720

Department of Mathematics, Vanderbilt University, Nashville, Tennessee 37240
*E-mail address*: `farley@msri.org`